\documentclass[12pt,reqno]{amsart}

\usepackage{amsmath,amsthm,amssymb,amsfonts, comment,fullpage}
\usepackage{braket}
\usepackage{mathtools}
\usepackage{bbm}
\usepackage{float}
\usepackage{caption}
\usepackage{times}
\usepackage{mathrsfs}
\usepackage{latexsym}
\usepackage[dvips]{graphics}
\usepackage{epsfig}
\usepackage{amsmath,amsfonts,amsthm,amssymb,amscd}
\input amssym.def
\input amssym.tex
\usepackage{color}
\usepackage{hyperref}
\usepackage{url}
\newcommand{\bburl}[1]{\textcolor{blue}{\url{#1}}}
\usepackage{diagbox}

\usepackage{tikz}
\usepackage{tkz-tab}
\usepackage{tkz-graph}
\usetikzlibrary{shapes.geometric,positioning}
\usetikzlibrary{decorations.pathreplacing}

\newcommand{\burl}[1]{\textcolor{blue}{\url{#1}}}

\numberwithin{equation}{section}

\newtheorem{thm}{Theorem}[section]

\newtheorem{defi}[thm]{Definition}

\theoremstyle{plain}

\newtheorem{definition}[thm]{Definition}

\newtheorem{theorem}[thm]{Theorem}

\newcommand\be{\begin{equation}}
\newcommand\ee{\end{equation}}
\newcommand\bee{\begin{equation*}}
\newcommand\eee{\end{equation*}}
\newcommand\bea{\begin{eqnarray}}
\newcommand\eea{\end{eqnarray}}
\newcommand\beae{\begin{eqnarray*}}
\newcommand\eeae{\end{eqnarray*}}
\newcommand\bi{\begin{itemize}}
\newcommand\ei{\end{itemize}}
\newcommand\ben{\begin{enumerate}}
\newcommand\een{\end{enumerate}}
\newcommand\bc{\begin{center}}
\newcommand\ec{\end{center}}
\newcommand\ba{\begin{array}}
\newcommand\ea{\end{array}}

\newcommand{\Z}{\ensuremath{\mathbb{Z}}}
\newcommand{\Q}{\mathbb{Q}}

\newcommand\frakfamily{\usefont{U}{yfrak}{m}{n}}
\DeclareTextFontCommand{\textfrak}{\frakfamily}

\newtheorem{rek}[thm]{Remark}

\newcommand{\hr}[1]{\href{#1}{\url{#1}}}

\newtheorem*{thm: 1 dimensional main result}{Theorem \ref{thm: 1 dimensional main theorem}}

\usepackage{placeins}

\usepackage{hyperref}
\usepackage{color}
\usepackage{url}

\title{Benford Behavior in Stick Fragmentation Problems}

\date{\today}

\author{Bruce Fang}
\email{\textcolor{blue}{\href{mailto:fangbaojun2002@gmail.com}{fangbaojun2002@gmail.com}}}
\address{Department of Mathematics, Williams College, Williamstown, MA 01267}

\author{Ava Irons}
\email{\textcolor{blue}{\href{mailto:ai4@williams.edu}{ai4@williams.edu}}}
\address{Department of Mathematics, Williams College, Williamstown, MA 01267}

\author{Ella Lippelman}
\email{\textcolor{blue}{\href{mailto:e_lippelman@coloradocollege.edu}{e\_lippelman@coloradocollege.edu}}}
\address{Department of Mathematics and Computer Science, Colorado College, Colorado Springs, CO 80903}

\author{Steven J. Miller}
\email{\textcolor{blue}{\href{mailto:sjm1@williams.edu}{sjm1@williams.edu}},  \textcolor{blue}{\href{Steven.Miller.MC.96@aya.yale.edu}{Steven.Miller.MC.96@aya.yale.edu}}}
\address{Department of Mathematics, Williams College, Williamstown, MA 01267}

\begin{document}

\pagenumbering{arabic}
\maketitle

Benford’s law is the statement that in many real-world data set, the probability of having digit \(d\) in base \(B\), where \(1\leq d\leq B\), as the first digit is \(\log_{B}\left((d+1)/d\right)\). We sometimes refer to this as weak Benford behavior, and we say that a data set exhibits strong Benford behavior in base \(B\) if the probability of having significand at most \(s\), where \(s\in[1,B)\), is \(\log_{B}(s)\). We examine Benford behaviors in stick fragmentation model. Building on the work on the 1-dimensional stick fragmentation model, we employ combinatorial identities on multinomial coefficients to reduce the high dimensional stick fragmentation model to the 1-dimensional model and provide a necessary and sufficient condition for the lengths of the stick fragments to converge to strong Benford behavior.

\tableofcontents

\section{Introduction}\label{sec: background}

\begin{figure}
    \centering
\includegraphics[width = 5in]{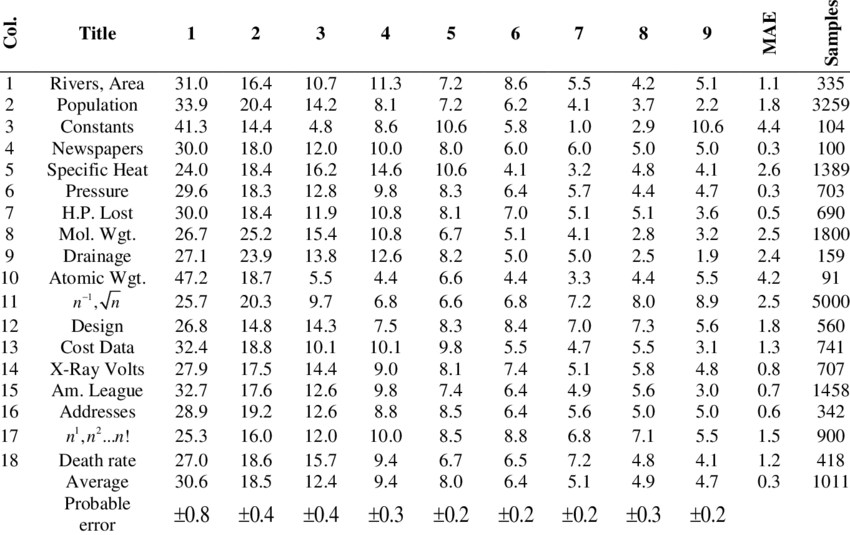}
\caption{From \cite{New}: Benford's $20,000$ observations.}
    \label{fig: benftablefig}
\end{figure}

In the late nineteenth century, astronomer Simon Newcomb observed that in the logarithmithic books at his workplace, certain pages were ``more worn than others" \cite{New}. In particular, there was more wear and tear in the earlier pages than the later pages. He deduced that there is a ``bias'' towards smaller leading digits, with the digit 1 showing up roughly $30\%$ of the time, the digit 2 showing up roughly $18\%$ of the time, and so on. Newcomb's findings were practically ignored until about fifty years later when physicist Frank Benford published his own research on the distribution of leading digits in \textit{Reader's Digest} \cite{Ben}. Benford displayed a table of roughly $20,000$ observations from twenty different sets of data, shown in Figure \ref{fig: benftablefig}. 

\begin{figure}
    \centering
\includegraphics[width = 5in]{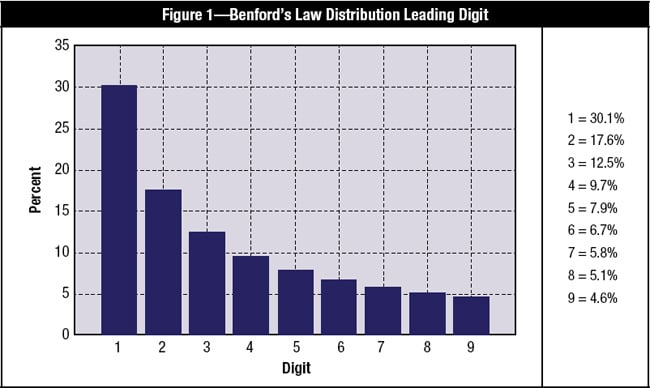}
\caption{The distribution of Benford's law.}
    \label{fig: benfordfigureintro}
\end{figure}

Benford called this distribution the ``law of Anomalous Numbers", but due to the popularity of his publication, the phenomena eventually became known as ``Benford's law"; see Figure \ref{fig: benfordfigureintro} for probabilities. Benford's law is a powerful phenomena that occurs in a variety of data, including accounting, elections, finance, geosciences, physics, population data, street addresses, and more. The law's prevalence makes it useful for ensuring data integrity, including using it as a method of fraud detection for tax returns, insurance claims, and expense reports \cite{Ni, PTTV}. For more information on the history of Benford's law, see \cite{BeHi, Hi1, Hi2, Mil, Rai}.

In 1986, Lemons \cite{Lemons} used Benford's law to analyze the partitioning of a conserved quantity. Since then, mathematicians and physicists have examined Benfordness of various fragmentation processes. Among these processes is the \emph{stick fragmentation process}, the subject of this paper. In the 1-dimensional multi-proportion stick fragmentation model, one starts with a stick of length $L$ and fixes an integer $m\geq 2$. At Stage 1, the stick is split into $m$ substicks of lengths $p^{(1, 1)}_1L, p^{(1, 1)}_2L, \dots, p^{(1, 1)}_mL$ according to proportions $p^{(1, 1)}_1, p^{(1, 1)}_2, \dots, p^{(1, 1)}_m$ drawn from a probability distribution on $(0,1)$, such that $p^{(1, 1)}_1+p^{(1, 1)}_2+\cdots+p^{(1, 1)}_m=1$. At Stage 2, each of the $m$ substicks created from Stage 1 is split into $m$ smaller substicks according to proportions $p_1^{(2, k)}, p_2^{(2, k)}, \dots, p_m^{(2, k)}$ drawn from the same probability distribution such that $p_1^{(2,k)}+p_2^{(2, k)}+\cdots+p_m^{(2, k)}=1$ for some $1\leq k\leq m$. In general, at Stage $N$, every substick created from Stage $N-1$ is further split into $m$ smaller sticks according to proportions $p^{(N, k)}_1, p^{(N, k)}_2, \dots, p^{(N, k)}_m$ drawn from the same probability distribution on $(0,1)$ such that $p^{(N, k)}_1+p^{(N, k)}_2+\cdots+p^{(N, k)}_m=1$ for some $1\leq k\leq m^{N-1}$. By the end of $N$ stages, there are $m^N$ sticks in total, and a question of particular interest to us is whether or not this fragmentation process results in stick lengths that are distributed according to Benford's law. In this paper, we specialize our investigation to a particular case of the multi-proportion stick fragmentation model, by fixing beforehand the $m$ proportions $p_1, p_2, \dots, p_m$ with which we split every substick. We provide a necessary and sufficient condition for such a fixed multi-proportion stick fragmentation model to result in sticks with lengths converging to a generalized Benford's law. 

\begin{theorem}\label{thm: 1 dimensional main theorem}
For any integer $m>2$, choose $p_1, p_2, \dots, p_{m-1}\in (0,1)$ such that $p_1+p_2+\dots+p_{m-1}<1$. Set $p_m:=1-(p_1+p_2+\dots+p_{m-1})$. At each stage, we cut a given stick according to proportions $p_1,p_2,\dots,p_{m-1}$ to create $m$ pieces. After $N$ iterations, we have $m^N$ sticks in total, of lengths
\begin{align}
A_{k_1,k_2,\dots,k_m} \ := \ Lp_1^{k_1}p_2^{k_2}\cdots p_m^{k_m},
\end{align}
for $0\leq k_1,k_2,\dots, k_m\leq N$ such that $k_1+k_2+\cdots+k_m=N$. Let $y_i=\log_{10}(p_i/p_{i+1})$ for $1\leq i\leq m-1$. Then the stick fragmentation model results in stick lengths that converge to strong Benford's law if and only if $y_i\not\in\Q$ for some $1\leq i\leq m-1$.
\end{theorem}

\subsection{Definitions and Theory of Benford's Law}
\label{sec: benftheory}

We begin with a definition of Benford's law (see for example \cite{Dia, Mil}).

\begin{defi}[Benford's Law for the Leading Digit] \label{defi: benflaw} We say a data set exhibits \textbf{Benford's law for the leading digit} if the frequency of leading digit $d$ is $\log_{10}(\frac{d+1}{d})$.
    
\end{defi}

There are many methods of proving that a data set follows Benford's law. A common one is to use the Uniform Characterization Theorem. An important definition for this characterization is the notion of the significand of a real number.

\begin{defi}[The Significand] \label{defi: sig} 
    For any positive $x$, we can express it as
    \be \label{eqn:sigeq}
    x \ = \ S_{10}(x)\cdot 10^{k_{10}(x)}
    \ee
    for unique $S_{10}(x)\in [1,10)$ and integer $k_{10}(x)$. We call $S_{10}(x)$ the \textbf{significand} of $x$.
\end{defi}

There is a more generalized version of Benford's law for the entire significand, not just for the first digit.
\begin{definition}
We say that a sequence of random variables $\{X^{(N)}\}_{N=1}^\infty$, converges to \textbf{strong Benford's law} if
\begin{equation}
\label{eq:StrongBenfordBehavior}
\lim_{N\rightarrow\infty}\mathbb{P}(S_{10}(X^{(N)}) \leq s) \ = \ \log_{10}(s)
\end{equation}
for all $s \in [1,10]$. Note that the increasing indices usually correspond to the growing size of the data set, so convergence to strong Benford's law is really an asymptotic statement.

\end{definition}
\begin{defi}[Uniform Distribution Modulo 1] \label{defi: unifmod1} A sequence of random variables $\{X^{(N)}\}_{N=1}^{\infty}$ converges to being equidistributed mod 1 if 
    \be \label{eqn:unifdis}
\lim_{N\rightarrow\infty}\mathbb{P}(X^{(N)}\textup{ mod }1\leq s) \ = \ s
    \ee
for all $s\in [0,1]$.
\end{defi}

Now, we are ready to state the Uniform Characterization Theorem \cite{Mil}.

\begin{thm}[Uniform Characterization Theorem] \label{thm:unithm} A sequence of random variables converges to strong Benford's law base 10 if and only if the sequence of the base 10 logarithms of the random variables converges to being equidistributed mod $1$.
\end{thm}

Thus, convergence to strong Benford's law and convergence to uniform distribution modulo 1 are in fact equivalent conditions. Thus, to prove (or disprove) convergence to strong Benford's law, it suffices to prove (or disprove) convergence to uniform distribution modulo 1.

\subsection{Fixed Proportion Stick Fragmentation Model and the Multinomial Distribution}
\label{sec: modelmulti}

Now that we have a proper foundation for Benford's law, it is time to describe in detail the model for this paper. The model can be thought of as an extension of a model from Becker et. al called the fixed single proportion stick fragmentation model \cite{Be}. Suppose we start with a stick of length $L$. We are going to split the stick at one fixed proportion $0<p<1$. After the first break, we have two sticks, namely of lengths $Lp$ and $L(1-p)$. Now we split the two sticks again at the same fixed proportion $p$, resulting in four sticks of lengths $Lp^2, Lp(1-p), L(1-p)p$, and $L(1-p)^2$. We continue this process for $N$ iterations. The fragmentation process follows a binomial distribution, with the  $2^N$ sticks of $N+1$ distinct lengths that follow a binomial distribution.

Becker et. al were interested in whether or not the leading digits of the significands of the stick lengths converge to strong Benford's law. They discovered a necessary and sufficient condition for the convergence to strong Benford's law and proved its necessity and sufficiency in \cite{Be}.

\begin{thm}[Fixed Single Proportion Stick Fragmentation Theorem, \cite{Be}]\label{thm:becker}
 Consider the fixed single proportion stick fragmentation model. Choose $y$ so that $10^y=(1-p)/p$. The fragmentation model results in stick lengths that converge to strong Benford's law if and only if $y\notin \Q$.
\end{thm}

For the non-Benford case, Becker et. al proved by noticing the cyclic behavior in the significands and by the multisection formula \cite{C}. By contrast, the Benford case was much harder to establish. They adopted methods from \cite{Dia}, \cite{KN}, and \cite{MT-B} to use truncation to show roughly equal probability between intervals, and to finally show equidistribution modulo 1. They were essentially able to prove that these stick lengths that followed a binomial distribution would converge to strong Benford's law if the ratio is equal to $10$ to an irrational power, and would not converge to Benford's law if the ratio is equal $10$ to a rational power. 

We extend their fixed single proportion stick fragmentation model to fixed multi-proportion stick fragmentation model. Becker et. al only explored the case when the stick $L$ is cut at one fixed proportion $p$ in every iteration. We push this case further, and see what happens when we cut the stick at multiple distinct fixed proportions $p_1,p_2,\dots,p_{m-1}$ in every iteration. Our stick model is as follows.

\begin{figure}
    \centering
\includegraphics[width = 6in]{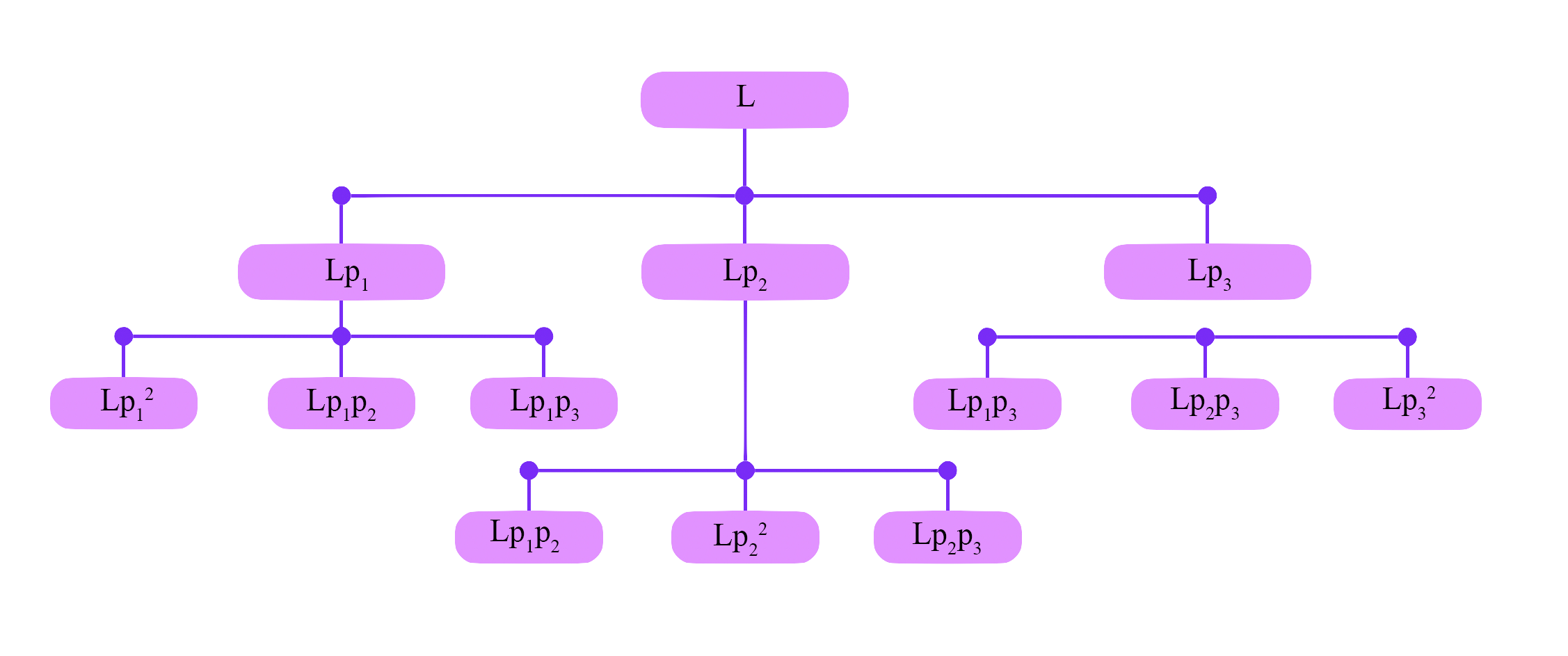}
\caption{A trinomial stick fragmentation with $m=3$ and $N=2$.}
    \label{fig: benffigtree}
\end{figure}

Suppose we have a stick of length $L$. We cut the stick simultaneously at fixed proportions $p_1,p_2,\dots,p_{m-1}\in (0,1)$, where $p_1+\cdots+p_{m-1}< 1$. Let $p_m:=1-(p_1+\cdots+p_{m-1})$.
Thus, after Stage 1, we are left with sticks of lengths $Lp_1,Lp_2,\dots,Lp_{m}$. At Stage 2, we cut each stick we gained from the previous iteration at the same fixed proportions $p_1,p_2,\dots,p_{m-1}$. Therefore, the sticks lengths that result from stage 2 are $Lp_1^2,Lp_1p_2,\dots Lp_1p_{m}$, $Lp_2p_1, Lp_2^2,\dots Lp_2p_{m}$, and so on. After Stage $N$, we are left with $m^N$ total sticks with $\binom{m+N-1}{N}$ distinct lengths (see Figure \ref{fig: benffigtree} for example when $m=3$ and $N=2$). We have this many distinct lengths as the problem can be thought as unordered sampling with replacement. Moreover, the stick lengths are distributed according to a generalized version of binomial distribution called the multinomial distribution, defined in terms of \textbf{multinomial coefficient}, which is a generalization of binomial coefficient. We introduce the following definition and result, which may be found in \cite{Ka}. 

\begin{defi}[Multinomial Coefficient] \label{defi: mcoef} For any non-negative integer $N$ and positive integer $m$, the multinomial coefficient is

\be \label{eq:mcoeq}
    \binom{N}{k_1,k_2,\dots,k_m} \ = \ \frac{N!}{k_1!k_2!\cdots k_m!}
\ee
for $k_1+k_2+\cdots+k_m=N$. A random vector $(X)=(X_1, X_2, \dots, X_m)$ follows a \textbf{multinomial distribution} with parameters $N$ and $(p):=(p_1, p_2, \dots, p_m)$ if
\begin{align}
\mathbb{P}((X)=(k)) \ = \ \binom{N}{k_1, k_2, \cdots, k_m}p_1^{k_1}p_2^{k_2}\cdots p_m^{k_m}
\end{align}
for all $(k):=(k_1, k_2, \dots, k_m)$ with $k_1+k_2+\cdots+k_m=N$.
\end{defi}

\begin{rek}
Formula \eqref{eq:mcoeq} represents the number of ways to choose $N$ objects with exactly $k_j$ objects of type $j$ when order does not matter.
\end{rek}

\begin{thm}[Multinomial Theorem] \label{thm: mthm} Let $N$ be any non-negative integer and $p_1,p_2,\dots,p_m$ be real numbers. Then
\be \label{eqn: meq}
(p_1+p_2+\cdots+p_m)^N \ = \ \:\sum_{k_1+k_2+\cdots+k_m\geq 0} \binom{N}{k_1, k_2,\dots,k_m}\:p_1^{k_1}p_2^{k_2}\cdots p_m^{k_m},
\ee
where the $k_j$'s are non-negative integers summing to $N$.
\end{thm}

\section{Fixed Multi-proportion Stick Fragmentation Model}
\label{sec:multithmsmall}

Recall that we are interested in studying whether or not a stick fragmentation process results in stick lengths that converge to strong Benford's law. For Becker et. al, they were able to prove that if the ratio $(1-p)/p$ is equal to $10$ to an irrational power, the stick lengths will follow strong Benford's law, but if the ratio is equal to $10$ to a rational power then the distribution of stick lengths will not follow strong Benford's law. 

In what follows, we generalize the results of Becker el. al to the fixed multi-proportion stick fragmentation model we introduced in Section \ref{sec: modelmulti}. In particular, we present Theorem \ref{thm: 1 dimensional main theorem} again along with the proof for the necessity of the condition.

\begin{thm: 1 dimensional main result}
For any integer $m>2$, choose $p_1, p_2, \dots, p_{m-1}\in (0,1)$ such that $p_1+p_2+\dots+p_{m-1}<1$. Set $p_m:=1-(p_1+p_2+\dots+p_{m-1})$. At each stage, we cut a given stick according to proportions $p_1,p_2,\dots,p_{m-1}$ to create $m$ pieces. After $N$ iterations, we have $m^N$ sticks in total, of lengths
\begin{align}
A_{k_1,k_2,\dots,k_m} \ := \ Lp_1^{k_1}p_2^{k_2}\cdots p_m^{k_m},
\end{align}
for $0\leq k_1,k_2,\dots, k_m\leq N$ such that $k_1+k_2+\cdots+k_m=N$. Let $y_i=\log_{10}(p_i/p_{i+1})$ for $1\leq i\leq m-1$. Then the stick fragmentation model results in stick lengths that converge to strong Benford's law if and only if $y_i\not\in\Q$ for some $1\leq i\leq m-1$.
\end{thm: 1 dimensional main result}

\begin{figure}[H] 
    \begin{minipage}{.5\textwidth}
        \centering
        \includegraphics[width=3 in]{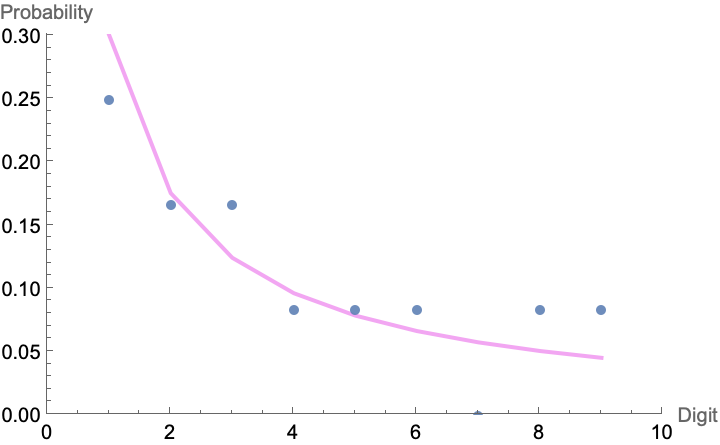}
        \captionof{figure}{\footnotesize $y_1=-1/2$,  $y_2=-1/3$, $y_3=-1/4$, and $N=100$.}
        \label{fig:quadsum1}
   \end{minipage}%
   \begin{minipage}{.5\textwidth}
          \centering
          \includegraphics[width= 3 in]{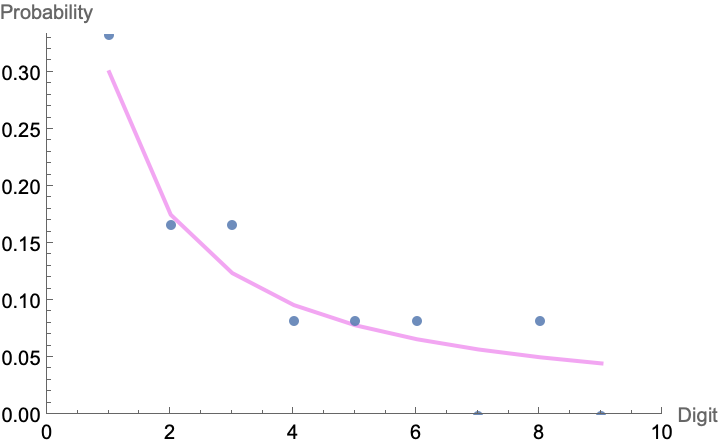}
          \captionof{figure}{\footnotesize $y_1=-1/4$, $y_2=-1/2$, $y_3=-1/6$, and $N=100$.}
          \label{fig:quadsum2}
   \end{minipage}
\end{figure}

\begin{figure}[H] 
    \begin{minipage}{.5\textwidth}
        \centering
        \includegraphics[width=3 in]{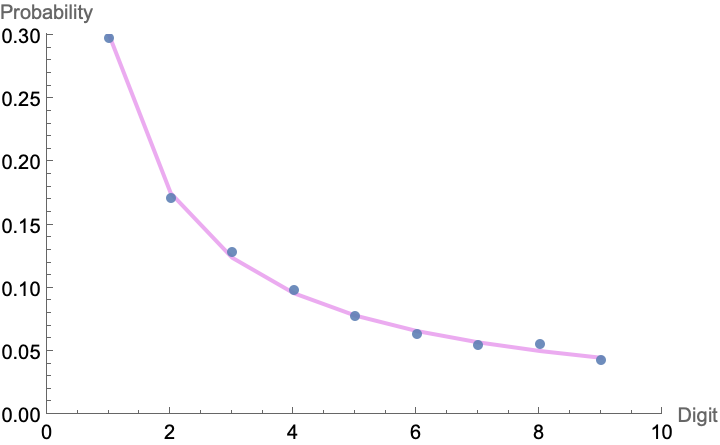}
        \captionof{figure}{\footnotesize $y_1=-1/2, \  y_2=-\sqrt{2},$  \\ and $N=1000$.}
        \label{fig:tri5}
   \end{minipage}%
   \begin{minipage}{.5\textwidth}
          \centering
          \includegraphics[width= 3 in]{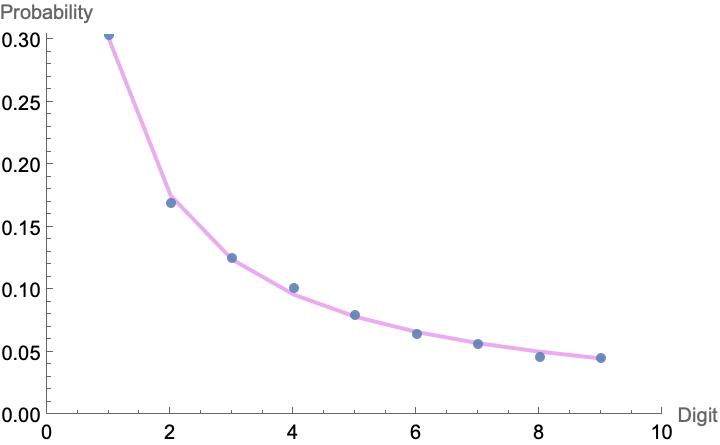}
          \captionof{figure}{\footnotesize $y_1=-1/3, \ y_2=-\sqrt{3}$, \\ and $N=1000$.}
          \label{fig:tri6}
    \end{minipage}
\end{figure}

\begin{figure}[H]
    \begin{minipage}{.5\textwidth}
        \centering        \includegraphics[width=3 in]{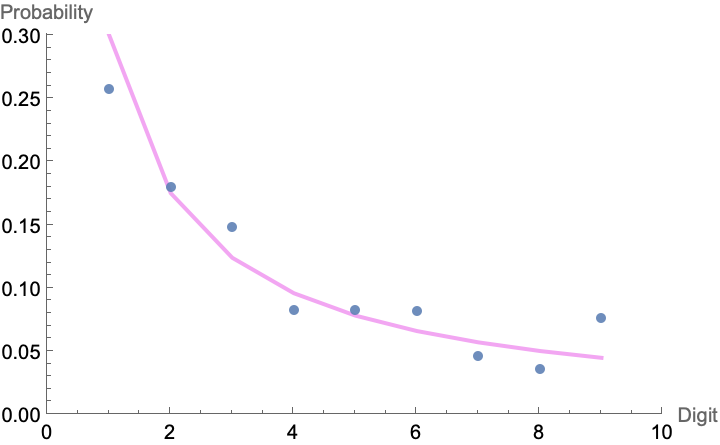}
        \captionof{figure}{\footnotesize $y_1=-\sqrt{2}$, $y_2=-1/3$, $y_3=-1/4$, and $N=100$.}
        \label{fig:quad1}
   \end{minipage}%
   \begin{minipage}{.5\textwidth}
          \centering  \includegraphics[width= 3 in]{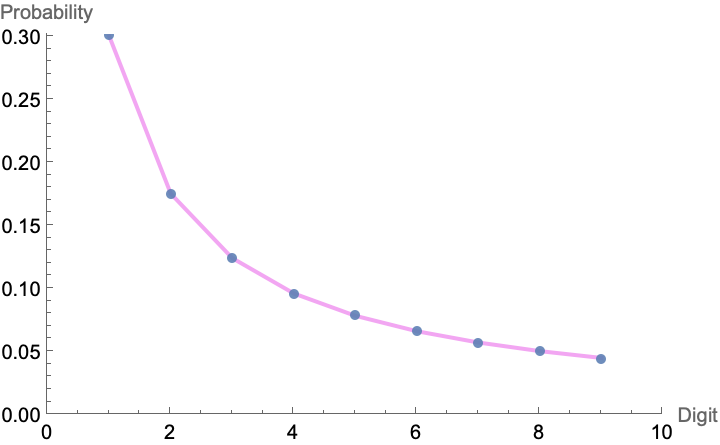}
          \captionof{figure}{\footnotesize $y_1=-\sqrt{3}$, $y_2=-1/10$, $y_3=-1/8$, and $N=100$.}
          \label{fig:quad2}
    \end{minipage}
\end{figure}

We prove the necessity of the condition, i.e., if $y_i\in \Q$ for all $1\leq i\leq m-1$, then the stick fragmentation model results in stick lengths that do not follow strong Benford's law. The proof of the sufficiency of the condition is long and technical, see \cite{FM} for the details. Before presenting the proof, we point the reader to some numerical simulation results in support of the theorem, as shown in the appendix (see Figures \ref{fig:tri1}, \ref{fig:tri2}, \ref{fig:quadsum1}, \ref{fig:quadsum2} for rational case and Figures \ref{fig:tri5}, \ref{fig:tri6}, \ref{fig:quad1}, \ref{fig:quad2} for irrational case).

\begin{proof}
To prove that the decomposition process results in stick lengths that do not converge to strong Benford's law, by the Uniform Characterization Theorem \ref{thm:unithm}, it suffices to show that
\begin{align}\log_{10}(A_{k_1, k_2, \dots, k_m}) \ = \ L\log_{10}\left(p_1^{k_1}p_2^{k_2}\dots p_m^{k_m}\right)
\end{align}
are not equidistributed mod 1. First, note that since Benford's law is scale invariant \cite{Hi1}, we can assume without loss of generality that $L=1$. We also notice that each stick length $A_{k_1,k_2,\dots, k_m}$ has the factorization
\begin{align}\label{factorization}
A_{k_1,k_2,\dots, k_m} \ = \ p_1^{k_1}p_2^{k_2}\cdots p_m^{k_m} \ = \ \left(\frac{p_1}{p_2}\right)^{k_1}\left(\frac{p_2}{p_3}\right)^{k_1+k_2}\cdots \left(\frac{p_{m-1}}{p_m}\right)^{\sum_{j=1}^{m-1}k_j}(p_m)^{N}.
\end{align}
Since $y_i\in \Q$ for all $1\leq i\leq m-1$, we can write $y_i=a_i/b_i$ for some $a_i\in\Z$ and $b_i\in \Z_{>0}$ and $\textup{gcd}(a_i, b_i)=1$, for all $1\leq i\leq m-1$. Thus
\begin{align}
A_{k_1,k_2,\dots, k_m} \ = \ \left(10^{\frac{a_1}{b_1}}\right)^{k_1}\left(10^{\frac{a_2}{b_2}}\right)^{k_1+k_2}\cdots \left(10^{\frac{a_{m-1}}{b_{m-1}}}\right)^{\sum_{j=1}^{m-1}k_j}(p_m)^N.
\end{align}
Taking the logarithm of both sides yields
\begin{align}
\log_{10}(A_{k_1,k_2,\dots, k_m}) \ = \ k_1\left(\frac{a_1}{b_1}\right)+(k_1+k_2)\left(\frac{a_2}{b_2}\right)+\cdots +\left(\sum_{j=1}^{m-1}k_j\right)\left(\frac{a_{m-1}}{b_{m-1}}\right)+N\log_{10}(p_m).
\end{align}
Note that the first term $k_1(a_1/b_1)$ on the RHS above has at most $b_1$ distinct values mod 1, since its values mod 1 are periodic with period at most $b_1$. More generally, for any $1\leq i\leq m-1$, the term $(\sum_{j=1}^{i}k_j)(a_i/b_i)$ has at most $b_i$ distinct values mod 1. Since $N\log_{10}(p_m)$ is constant for all $A_{k_1,k_2,\dots, k_m}$, then $\log_{10}(A_{k_1,k_2,\dots, k_m})$ has at most $\prod_{i=1}^{m-1}b_i$ distinct values. Moreover, this number does not depend on $N$ and is finite for fixed $m$, $b_1, \dots, b_{m-1}$. Since a uniform distribution is continuous and can take any values in $[0,1]$, which is uncountable, then $\log_{10}(A_{k_1, k_2, \dots, k_m})$ is not equidistribted mod 1. Thus, the decomposition process results in stick lengths that do not converge in distribution to strong Benford's law.
\end{proof}

\section{Future Work}
While our proof for the necessity of the condition of Theorem \ref{thm: 1 dimensional main theorem} shows that $\log_{10}(A_{k_1, k_2, \dots, k_m})$ has at most $\prod_{i=1}^{m-1}b_i$ distinct values and thus does not converge to being equidistributed mod 1, it does not explicitly determine the distribution of $\log_{10}(A_{k_1, k_2, \dots, k_m})$. It would be an interesting question to see if the distribution converges to discrete uniform and if not, to characterize the distribution under various general assumptions.

\section*{Acknowledgments}
This work was partially supported by Williams College Summer Science Program Research Fellowship, the Finnerty Fund, and NSF Grant DMS2241623. The third author thanks her thesis advisor Molly Moran for guidance and encouragement.


\newpage

\section{Numerical Evidence for Theorem \ref{thm: 1 dimensional main theorem}}

\subsection{Evidence for Non-Benford Behavior} The following figures provide numerical evidence for non-Benford behavior when the exponents are all rational.

\begin{figure}[H]

    \begin{minipage}{.5\textwidth}
        \centering
        \includegraphics[width=3 in]{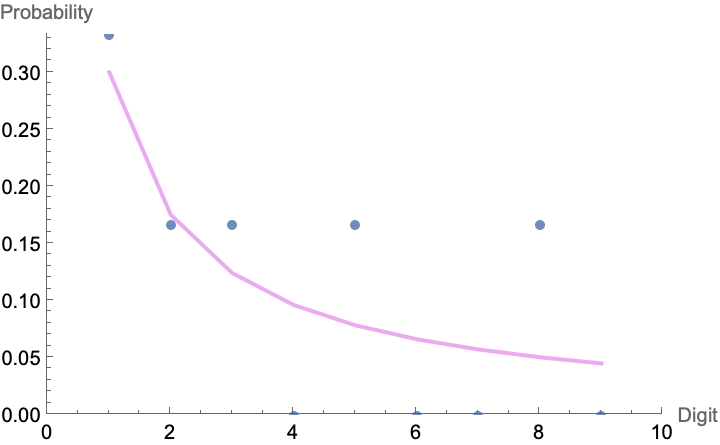}
        \captionof{figure}{\footnotesize $y_1=-1/3, \  y_2=-1/2,$  \\ and $N=1000$.}
        \label{fig:tri1}
   \end{minipage}%
   \begin{minipage}{.5\textwidth}
          \centering
          \includegraphics[width= 3 in]{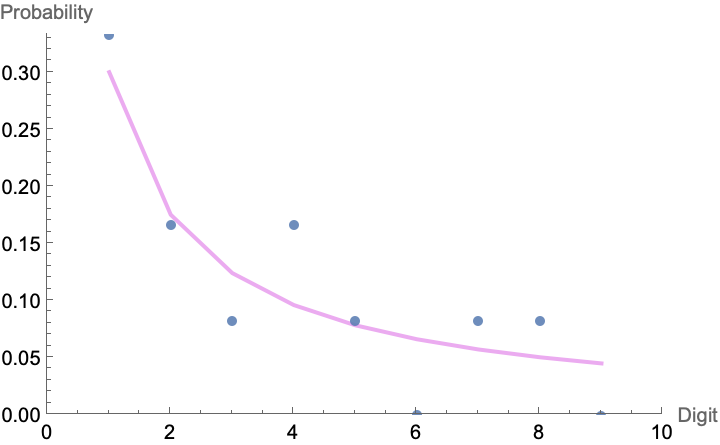}
          \captionof{figure}{\footnotesize $y_1=-1/4, \ y_2=-1/6,$ \\ and $N=1000$.}
          \label{fig:tri2}
    \end{minipage}
\end{figure}

\subsection{Evidence for Benford Behavior} The following figures provide numerical evidence for Benford behavior when at least one exponent is irrational.

\end{document}